\documentclass[a4paper]{amsart}
\usepackage{amsmath,amsthm,amssymb}
\usepackage{hyperref}
\usepackage{fullpage}

\theoremstyle{plain}
\newtheorem*{mainthm}{Main Theorem}
\theoremstyle{remark}
\newtheorem*{rmk}{Remark}

\newcommand{\bC}{\mathbb{C}}
\newcommand{\cC}{\mathcal{C}}
\newcommand{\bF}{\mathbb{F}}
\newcommand{\cG}{\mathcal{G}}
\newcommand{\bZ}{\mathbb{Z}}
\DeclareMathOperator{\Rep}{Rep}
\DeclareMathOperator{\Irr}{Irr}
\newcommand{\Path}{\mathrm{Path}}

\begin{document}
	\title{A note on quantum subgroups of free quantum groups}
	\author{Mao Hoshino and Kan Kitamura}
	\address{Graduate School of Mathematical Sciences, the University of Tokyo, 3-8-1 Komaba Meguro-ku Tokyo 153-8914, Japan}
	\email{MH: mhoshino@ms.u-tokyo.ac.jp \qquad
		KK: kankita@ms.u-tokyo.ac.jp}
    \subjclass{Primary 46L67, Secondary 05C25}
    \thanks{MH was supported by Forefront Physics and Mathematics Program and JSPS KAKENHI Grant Number JP23KJ0695. KK was supported by JSPS KAKENHI Grant Number JP22KJ0618.}
	\begin{abstract}
		In this short note, quantum subgroups in finite free products of the Pontryagin duals of free unitary quantum groups are classified. They correspond to pairs of a subgroup $\Gamma$ and a subset $S$ of the free group $\mathbb{F}_n$ such that $S$ is $\Gamma$-invariant, containing $\Gamma$, and connected in the Cayley graph of $\mathbb{F}_n$. 
	\end{abstract}
	\maketitle

	\section{Introduction}
	
	In this note, we classify quantum subgroups in the Pontryagin duals of free unitary quantum groups or, more generally, in their finite free products. 
	Associated with $Q\in GL_N(\bC)$ for $N\geq 2$, van Daele--Wang~\cite{vanDaele-Wang:universal} defined a compact quantum group called a free unitary quantum group, which shall be denoted by $U^+_Q$. 
	It is constructed by forgetting the commutativity relations for the coefficients of a unitary matrix. Variations of this construction 
	have provided abundant examples of compact quantum groups. 
	
	For general conventions and fundamental facts on compact quantum groups and tensor categories, we refer \cite{Neshveyev-Tuset:book}. 
	For a compact quantum group $G$, its finite dimensional unitary representations form a rigid semisimple C*-tensor category with a simple unit object, denoted by $\Rep G$. We write $\Irr G$ for the set of all isomorphism classes of simple objects in $\Rep G$. 
	We write $\widehat{G}$ for the Pontryagin dual of $G$, which is a discrete quantum group. 
	
	The notion of free product for (the C*-algebras of) discrete quantum groups was studied by Wang~\cite{Wang:free}. 
	The main object of this note is the free product ${\ast}_{j=1}^n\widehat{U}^+_{Q_j}$ for $n\in\bZ_{\geq 1}$, $Q_j\in GL(N_j,\bC)$, $N_j\geq 2$, defined by the Pontryagin dual $\widehat{F}$ of the compact quantum group $F$ given by $C_u(F)=\ast_{j=1}^nC_u(U_{Q_j}^+)$ with a canonical comultiplication. 
	This quantum group $\widehat{F}$ has the following remarkable property analogous to the free group $\bF_n$ with $n$ generators: for any compact quantum group $G$, any $n$-tuple of unitary representations $(u_j)_{j=1}^{n}$ of $G$ with $\dim_\bC u_j = N_j$ and $Q_j^{*-1}u_j^{\top}Q_j^*$ being a unitary uniquely induces a unital $*$-homomorphism $C_u(F)\to C_u(G)$ preserving comultiplications. 
	Quantum subgroups of $\widehat{F}$ has been recently classified by Freslon--Weber~\cite{Freslon-Weber:discrete} when $\widehat{F}=\widehat{U}^+_Q$, together with rich structural results for such quantum subgroups. 
	We classify quantum subgroups of $\widehat{F}$ in general, based on a simple observation relating $\Irr F$ to the path space of the Cayley graph of $\bF_n$. 
	Our classification exhibits to what extent more quantum subgroups exist than mere subgroups of free groups and thus could be regarded as a complementary result to \cite{Freslon-Weber:discrete} towards the quantum analogue of Kurosh's theorem.

	\section{Classification}
	Note that the classification of quantum subgroups in $\widehat{F}$ is equivalent to the classification of idempotent complete full rigid C*-tensor subcategories in $\Rep F$, as a consequence of Woronowicz's Tannaka--Krein duality. Since the latter is a purely categorical problem and in fact relies only on the fusion rule of $\Rep F$, we work with $\Rep F$ rather than $\widehat{F}$ and assume for simplicity that $\Rep F$ is skeletal by strictification. 
	
	The theory of unitary representations of $U^+_Q$ was investigated by Banica~\cite{Banica:groupe}. We recall its fusion rule. 
	Consider the free product of additive monoids $M:=\bZ_{\geq 0}\ast \bZ_{\leq 0}$. 
	When we regard $r\in \bZ_{\geq 0}$ and $s\in\bZ_{\leq 0}$ as elements in $M$, we write $[r]$ and $[s]$. 
	Then, we can identify $\Irr U^+_Q$ with $M$ so that 
	the unit object is $[0]$, $\overline{[r]}=[-r]$ for $r \in \bZ$, $\overline{(xy)}=\bar{y}\bar{x}$ for $x,y \in M$, and the tensor product is determined recursively by
	\begin{align*}
		&
		(x[r])\otimes([s]y) = \left\{ \begin{array}{ll}
			x[r][s]y=x[r+s]y & \text{if $rs=1$,} \\
			x[r][s]y \oplus ( x\otimes y )& 
			\text{if $rs=-1$}
		\end{array} \right. 
	\quad (x,y \in M, r,s \in \{\pm 1\}). 
	\end{align*}
	Combined with \cite{Wang:free}, we can canonically identify $\Irr F$ with the monoid $M^{\ast n}$ as a set. 
	For $1\leq j\leq n$, we put $\alpha_j\in\bZ_{\geq 0}$, $\beta_j\in\bZ_{\leq 0}$ for the generators of the $j$th component $M$ in $M^{\ast n}$, and $a_j$ for the generator of the $j$th component $\bZ$ in $\bF_n$. 
	We shall write $[a_j^r]:=\alpha_j^r$ and $[a_j^{-r}]:=\beta_j^r$ in $M^{\ast n}$ for $r\in \bZ_{\geq 0}$.

	Consider the Cayley graph $\cG$ of $\bF_n$, where the set of its oriented edges is $\{ (g,ga_j^{\pm 1}) \mid g\in \bF_n, 1\leq j\leq n \}$. Let $\Path(\cG)$ be the set of all finite paths in $\cG$ starting from $1_{\bF_n}$, possibly with turning-backs. 
	We can identify $\Path(\cG)$ with $M^{\ast n}$ by assigning each path of the form $(g_0,g_1,\cdots,g_k)\in \bF_n^{k+1}$ with $k\geq 0$, $g_0=1_{\bF_n}$, $g_l=g_{l-1}a_{j_l}^{r_l}$, $1\leq j_l\leq n$, $r_l=\pm1$ for $1\leq l\leq k$ to $[a_{j_1}^{r_1}][a_{j_2}^{r_2}]\cdots[a_{j_k}^{r_k}]\in M^{\ast n}$. 
	Then, extending our convention of $[a_j^r]$ for $r\in\bZ$, we define $[g]\in M^{\ast n}$ for $g\in \bF_n$ 
	by the shortest path from $1_{\bF_n}$ to $g$ via the identification $M^{\ast n}=\Path(\cG)$. 
	From now, we promise that the symbol $k$, $j$ and $r$, possibly accompanied with indices, stand for elements in $\bZ_{\geq 0}$, $\bZ\cap[1,n]$, and $\{\pm1\}$ respectively, unless clarified otherwise. 
	Note that any element in $\Irr F = M^{\ast n}$ can be expressed in the form $[a_{j_1}^{r_1}][a_{j_2}^{r_2}]\cdots [a_{j_k}^{r_k}]$ without redundancy.
	
	\begin{rmk}
		We reinterpret the fusion rule of $\Rep F$ in terms of $\cG$. 
		Take any $e,e'\in\Path(\cG)$, whose endpoints are denoted by $g,g'\in\bF_n$, respectively. 
		We write $ge'$ for the path from $g$ to $gg'$ given by the left translation of $\bF_n$ on $\cG$ and $ee'\in\Path(\cG)$ for the concatenation of $e$ and $ge'$. 
		Note that $ee'$ coincides with the product of the monoid $M^{\ast n}$. 
		We set $e_0:=e$, $e'_0:=e'$ and define $e_{l+1},e'_{l+1}\in\Path(\cG)=M^{\ast n}$ recursively on $l\in\bZ_{\geq 0}$ so that $e_{l+1}[a_j^{r}]=e_{l}$ and $[a_j^{-r}] e'_{l+1}= e'_{l}$ as long as such $j$ and $r$ exist. If we regard $e,e'\in\Irr F$, the tensor product $e\otimes e'$ is the direct sum of all paths of the form $e_l e_l'\in\Path(\cG)=\Irr F$. 
		Also, when we write $\overline{e}\in\Path(\cG)$ for the path from $1_{\bF_n}$ to $g^{-1}$ given by reversing $e$, the involution $e\mapsto\overline{e}$ corresponds to the conjugate operation on $\Irr F$. 
		If we put $V(e)\subset \bF_n$ for the set of all vertices appearing in $e$, we have $V(e_l e'_l)\subset V(e e')=V(e)\cup gV(e')$ for $e_l$ and $e_l'$ above, 
		and $V(\overline{e})=g^{-1}V(e)$. 
	\end{rmk}
	
	Our main result is as follows. 
	
	\begin{mainthm}\renewcommand{\theequation}{\Roman{equation}}
		For any pair $(\Gamma,S)$ of a subgroup $\Gamma\subset \bF_n$ and a left $\Gamma$-invariant subset $S\subset \bF_n$ containing $\Gamma$ and connected in the Cayley graph $\cG$, 
		there is an idempotent complete full rigid C*-tensor subcategory of $\Rep F$ such that the set of its irreducible objects equals 
		\begin{align}\label{eqmain}
			\bigl\{ [a_{j_1}^{r_1}][a_{j_2}^{r_2}]\cdots [a_{j_k}^{r_k}]\in\Irr F 
			\,\big|\, 
			a_{j_1}^{r_1}a_{j_2}^{r_2}\cdots a_{j_k}^{r_k}\in \Gamma, 
			a_{j_1}^{r_1}a_{j_2}^{r_2}\cdots a_{j_l}^{r_l}\in S 
			\text{ for all }0\leq l\leq k 
			\bigr\}. 
		\end{align}
		
		Conversely, any idempotent complete full rigid C*-tensor subcategory of $\Rep F$ is of this form for a pair $(\Gamma,S)$ with the conditions above, which must be necessarily unique. 
	\end{mainthm}\setcounter{equation}{0}
	
	\begin{proof}
		It is not hard to see the first half of the statement with the aid of the remark above. 
		We show the latter half. 
		Let $\cC$ be an idempotent complete full rigid C*-tensor subcategory of $\Rep F$ and $I$ be the set of irreducible objects in $\cC$. 
		For any $[a_{j_1}^{r_1}]\cdots [a_{j_k}^{r_k}]\in I \setminus \{[1_{\bF_n}]\}$, 
		inductively on $l=1,\cdots,k$, we see 
		\begin{align}
			&\label{eq1}
			[a_{j_1}^{r_1}\cdots a_{j_l}^{r_l}][a_{j_{l+1}}^{r_{l+1}}]\cdots [a_{j_k}^{r_k}]
			\leq 
			\bigl( [a_{j_1}^{r_1}\cdots a_{j_{l-1}}^{r_{l-1}}] [a_{j_{l-1}}^{-r_{l-1}}\cdots a_{j_1}^{-r_1}] \bigr)
			\otimes 
			\bigl( [a_{j_1}^{r_1}\cdots a_{j_{l-1}}^{r_{l-1}}] [a_{j_l}^{r_l}]\cdots [a_{j_k}^{r_k}] \bigr)
			\in\cC,\text{ and} 
			\\&\label{eq2}
			[a_{j_1}^{r_1}\cdots a_{j_{l}}^{r_{l}}] [a_{j_{l}}^{-r_{l}}\cdots a_{j_1}^{-r_1}]
			\leq 
			\bigl( [a_{j_1}^{r_1}\cdots a_{j_l}^{r_l}][a_{j_{l+1}}^{r_{l+1}}]\cdots [a_{j_k}^{r_k}] \bigr) 
			\otimes 
			\bigl( [a_{j_k}^{-r_k}]\cdots [a_{j_{l+1}}^{-r_{l+1}}][a_{j_l}^{-r_l}\cdots a_{j_1}^{-r_1}] \bigr) 
			\in\cC. 
		\end{align}
		
		Clearly, $S:=\{ a_{j_1}^{r_1}\cdots a_{j_l}^{r_l} \mid [a_{j_1}^{r_1}]\cdots [a_{j_k}^{r_k}]\in I, 0\leq l\leq k \} 
		\subset \bF_n$ is connected in $\cG$, 
		and 
		$\Gamma := \{ a_{j_1}^{r_1}\cdots a_{j_k}^{r_k} \mid [a_{j_1}^{r_1}]\cdots [a_{j_k}^{r_k}]\in I \}$ is a subgroup of $\bF_n$ contained in $S$. 
		Then, (\ref{eq1}) and (\ref{eq2}) for $l=k$ imply $[g]\in I$ for all $g\in \Gamma$ and $[g][g^{-1}]\in I$ for all $g\in S$, respectively. 
		For all $g\in \Gamma$ and $h\in S$, it holds $[g][h][h^{-1}]\in I$ and thus $gh\in S$. 
		By construction, the set defined by (\ref{eqmain}) contains $I$. They coincide because 
		for any $[a_{j_1}^{r_1}]\cdots [a_{j_k}^{r_k}]\in \Irr F$ satisfying the conditions in (\ref{eqmain}), 
		\begin{align*}
			&
			[a_{j_1}^{r_1}]\cdots [a_{j_k}^{r_k}] 
			\leq 
			\bigl( [a_{j_1}^{r_1}]\cdots [a_{j_k}^{r_k}][a_{j_k}^{-r_k}\cdots a_{j_1}^{-r_1}] \bigr) 
			\otimes 
			[a_{j_1}^{r_1} \cdots a_{j_k}^{r_k}]
			\\\leq{}&
			\bigl([a_{j_1}^{r_1}] [a_{j_1}^{-r_1}]\bigr) 
			\otimes 
			\bigl([a_{j_1}^{r_1}a_{j_2}^{r_2}] [a_{j_2}^{-r_2}a_{j_1}^{-r_1}]\bigr) 
			\otimes 
			\cdots
			\otimes 
			\bigl([a_{j_1}^{r_1}\cdots a_{j_k}^{r_k}] [a_{j_k}^{-r_k}\cdots a_{j_1}^{-r_1}]\bigr) 
			\otimes 
			[a_{j_1}^{r_1} \cdots a_{j_k}^{r_k}] 
			\in\cC. 
		\qedhere\end{align*}
	\end{proof}
	
\subsection*{Acknowledgments}
	The authors would like to thank Amaury Freslon and Moritz Weber for their comments on the draft and their encouragement to release the content. They would like to thank their advisor Yasuyuki Kawahigashi for his invaluable support. 

\end{document}